\documentclass[11 pt]{amsart}
\usepackage{amssymb}
\usepackage{amsmath}
\usepackage{amsfonts}
\usepackage{graphicx}
\usepackage{amsthm}
\usepackage{enumerate}
\usepackage[mathscr]{eucal}
\usepackage{verbatim}
 \theoremstyle{plain}

\theoremstyle{plain}

\theoremstyle{remark}

\def\bbR{{\mathbb {R}}}

\def\cC{{\mathcal {C}}}

\def\cT{{\mathcal {T}}}

\def\cL{{\mathcal {L}}}

\begin{document}
  
\date{June, 2007}  

\title
{Convolution estimates and model surfaces of low codimension}

\author[]
{Daniel M. Oberlin}

\address
{Department of Mathematics \\ Florida State University \\
Tallahassee, FL 32306}
\email{oberlin@math.fsu.edu}


\thanks{The author was supported in part by NSF grant DMS-0552041.}

\begin{abstract}
For $k\geq d/2$ we give examples of measures on $k$-surfaces in $\bbR
^d$. These measures satisfy convolution estimates which are nearly
optimal. \end{abstract}

\maketitle


Suppose that $S$ is a smooth $k$-dimensional surface in $\bbR ^d$ and that
$\mu$ is a smooth positive Borel measure on $S$. Suppose further that
$\mu$ satisfies the convolution estimate 
\begin{equation}\label{LpLq}
\|\mu\ast f\|_{L^q (\bbR ^d )}\leq C\, \|f\|_{L^p (\bbR ^d )},
\end{equation}
where the norms are computed using Lebesgue measure $m_d$ on $\bbR ^d$.
Then it is well-known that $p\leq q$. Convolution with the characteristic 
function of a small ball shows that $(1/p,1/q)$ must lie in the triangle 
$\Delta (k,d)$ with vertices $(0,0)$, $(1,1)$, and 
$\big(d/(2d-k),(d-k)/(2d-k)\big)$. And a result of Ricci (\cite{R}), which
extends an observation of Carbery and Christ, shows that if $k(k+3)<2d$,
then \eqref{LpLq} also implies that 
\begin{equation}\label{Ricci}
\frac{1}{p}-\frac{1}{q}\leq \frac{2k}{6d-k^2-5k}.
\end{equation}
Let $\cT (k,d)$ be $\Delta (k,d)$ if $k(k+3)\geq 2d$ and the subset of 
$\Delta (k,d)$ defined by \eqref{Ricci} if $k(k+3)< 2d$. Suppose now that 
$S$ has the form 
\begin{equation}\label{param}
\{\big(y ;\Phi _1 (y), \Phi _2 (y),\dots ,\Phi _l (y)\big):
y\in G\}
\end{equation}
where $G$ is a nonempty open subset of $\bbR
^k$, where $l=d-k$, and where the functions $\Phi _j :\bbR ^k \rightarrow
\bbR$ are homogeneous polynomials. Let $\mu$ be the measure on $S$ induced
by $m_k$ on $G$. Then we will say that $S$ is a {\it model
surface} if \eqref{LpLq} holds whenever $(1/p,1/q)$ lies in the interior
of $\cT (k,d)$. 

\noindent{\bf Examples:}

(i) the paraboloids $\{(y;|y|^2 ):y\in \bbR ^{d-1}, |y|<1\}$ (see, e.g., 
pp. 370--371 in \cite{Stein});

(ii) the moment curves $\{(y;y^2 ,\dots ,y^d ):0< y<1\}$ (see 
\cite{Christ});

(iii) the monomial surfaces $\{\big(y;\Phi _1 (y),\dots ,\Phi _l (y) \big)
: y\in \bbR ^k ,|y|<1\}$ where $l=k+\frac{k(k-1)}{2}$ and the functions 
$\Phi _j$ are the distinct quadratic monomials (see \cite{R});

(iv) the $3$-surface $\{ (y_1 ,y_2 ,y_3 ;y_1^2 + y_2^2 ,y_2^2 + y_3^2 ):0<
y_j <1\}$ 
in $\bbR ^5$ (see \cite{O});

(v) certain surfaces of the form $\{\big(y;\Phi _1 (y),\dots ,\Phi _l
(y) \big) : y\in \bbR ^k ,|y|<1\}$ where $l=k$ (see \cite{DG}).

\noindent Of course most polynomial surfaces $S$
of the form \eqref{param} are not model surfaces in our sense: the
convolution requirement rules out degeneracies which result from the 
presence of \lq\lq flatness" or the lack of \lq\lq curvature". When
$k=1$ or $k=d-1$ there are obvious and simple technical interpretations of
\lq\lq curvature". In a few other cases there are technical
interpretations which are neither obvious nor simple. 
For example, when $k=2$ and $d=4$ the interpretation is that 
\begin{equation*}\begin{aligned}
(\Phi_{1y_1 y_1}\Phi_{2y_1 y_2}-\Phi_{2y_1 y_1}
\Phi_{1y_1  y_2})
&
(\Phi_{2y_2 y_2}\Phi_{1y_1 y_2}-\Phi_{2y_1 y_2}\Phi_{1y_2 y_2})
\\
&
-\big((\Phi_{1y_1 y_1}\Phi_{2y_2 y_2}-\Phi_{2y_1 y_1}\Phi_{1y_2
y_2})\big)^2
\end{aligned}\end{equation*}
not vanish. At any rate, the examples mentioned above, along with certain
of their Cartesian products, constitute a fairly 
complete list of the known model surfaces. The aim of this note is to
extend that list by providing examples of model surfaces whenever 
$k\geq \frac{d}{2}$. 

Fix positive integers
$k$ and $l$ with $1\leq l\leq k$ and put $d=k+l$. Let $\cC =[c_i^j ]$ be
a $k$ by
$l$ matrix of real numbers. For $1\leq j\leq l$ define bilinear forms 
$L_j : \bbR ^k \times \bbR ^k \rightarrow \bbR$ by 
$$
L_j (x,y)=\sum_{i=1}^k c_i^j x_i y_i
$$
and put $\Phi _j (y) =L_j (y,y)$. We will say that $\cC$ satisfies
condition (*) if every $l$ by $l$ submatrix of $\cC$ is nonsingular.

\noindent{\bf Theorem.} {\it
With the $\Phi _j$ as above, with $G=B(0,1)$, and with
$S$ given by \eqref{param}, suppose that (*) holds. Then $S$ is a model
surface.
}

\noindent{\it Proof of theorem:} Since $k(k+3)\geq 2d$, it is required to
establish \eqref{LpLq} whenever 
$(\frac{1}{p},\frac{1}{q})$ lies in the interior of $\Delta (k,d)$. With
$q_0 =\frac{2d-k}{d-k}$, an interpolation argument shows that it is
suffices to prove that
\begin{equation}\label{convineq}
\|\mu\ast \chi _E \|_{q_0} \leq C({p})\, m_d
(E)^{1/{p}}
\end{equation} 
for measurable $E\subset \bbR ^d$ and ${p}>\frac{2d-k}{d}$. And,
since $\mu$ has compact support, we
can also assume that $E\subset B(0,1)$. For such $E$, \eqref{convineq}
will follow, as in \cite{O}, from the auxiliary inequality
\eqref{ineq2} below. Thus, writing $\Phi (y)=\big(y;\Phi _1 (y),\dots
,\Phi _l (y)\big)$, 
\begin{equation}\label{ineq1}
\begin{aligned}
&
\|\mu\ast\chi _E \|^{q_0}_{q_0}
\\
&
=\int_{\bbR ^d}\int_{B(0,1)} \chi _E \big(z-\Phi (x)\big)\,dm_k (x)
\Big(\int_{B(0,1)} \chi _E \big(z-\Phi (y)\big)\,dm_k (y)\Big)^{q_0
-1}dm_d (z)
\\
&
=\int_{\bbR ^d} \chi _E (z)\int_{B(0,1)}\Big(\int_{B(0,1)} \chi _E
\big(z+\Phi (x)-\Phi (y)\big)\,dm_k
(y)\Big)^{q_0 -1}dm_k (x)\,dm_d (z).
\end{aligned}\end{equation}
Now assume, for the moment, the inequality 
\begin{equation}\label{ineq2}
\Big( \int_{B(0,1)} \Big[ \int_{B(0,1)} \chi _{\widetilde{E}}
\big( \Phi (x)-\Phi (y)\big)\,dm_k
(y)\Big]^{d/l}
dm_k (x)\Big)^{l/d}\leq C(\widetilde{p} )\,m_d
(\widetilde{E})^{1/{\widetilde{p}}}
\end{equation}
for $\widetilde{p} >\frac{d}{k}$ and $\widetilde{E}\subset B(0,2)$.
Since $q_0 -1 =\frac{d}{l}$, \eqref{ineq1} and \eqref{ineq2} yield
\begin{equation*}
\|\mu\ast\chi _E \|_{q_0}\leq C(\widetilde{p} )\,
m_d(E)^{(1+\frac{1}{\widetilde{p}}\frac{d}{l})\frac{1}{q_0}}. 
\end{equation*}
If $\frac{1}{{p}}=(1+\frac{1}{\widetilde{p}}\frac{d}{l})\frac{1}{q_0}$,
then $\widetilde{p} >\frac{d}{k}$ if and only if
${p}>\frac{2d-k}{d}$. Thus, as claimed,  
\eqref{convineq} will follow from \eqref{ineq2}. Now \eqref{ineq2}
is equivalent to the inequality, for nonnegative $f$,
\begin{equation}\label{ineq3}\begin{aligned}
\int_{B(0,1)} \int_{B(0,1)} f(x)\, \chi
_{\widetilde{E}}\big(x-y;\sum\nolimits c_i^1
(x_i^2-y_i^2 ),\dots ,\sum\nolimits 
&
c_i^l (x_i^2 
-y_i^2 )\big)\, 
dm_k (y) \, dm_k (x) 
\\
&
\leq C(\widetilde{p} )\,\|f\|_{L^{d/k}(\bbR ^k )}m_d (\widetilde
{E})^{1/\widetilde{p} },
\end{aligned}\end{equation}
where $\sum$ means $\sum_{i=1}^k$. In the $y$-integral
we change variables to obtain 
\begin{equation*}
\int_{B(0,1)}\int_{B(0,1)}
f(x)\, \chi _{\widetilde{E}}\big(y ;\sum\nolimits c_i^1
(2x_i y_i -y_i^2),\dots
,\sum\nolimits c_i^l (2x_i y_i -y_i^2)\big)\, dm_k (y) \, dm_k (x).
\end{equation*}
If $E\subset \bbR ^d$ is defined by 
\begin{equation*}\label{ineq4}
\chi _E (y_1 ,\dots ,y_k ; u_1 ,\dots ,u_l)=
\chi _{\widetilde{E}} \big(y_1 ,\dots ,y_k , 2u_1 - \sum c_i^1 y_i^2,
\dots ,2u_l - \sum c_i^l y_i^2 \big),
\end{equation*}
then $m_d (E)=2^{-l}m_d (\widetilde {E})$ 
and the left hand side of \eqref{ineq3} may be written 
$$
\int_{B(0,1)}
\int_{B(0,1)} f(x) \chi _E \big(y; L_1 (x,y),\dots ,L_l (x,y)\big)\,dm_k
(y)\,dm_k
(x).
$$
Thus \eqref{ineq3} will follow from 
\begin{equation}\label{ineq5}
\int_{B(0,1)}
\int_{B(0,1)} f(x) \chi _E \big(y; L_1 (x,y),\dots ,L_l (x,y)\big)\,dm_k
(y)\,dm_k
(x)\leq C
\,\|f\|_{L^{d/k}(\bbR ^k )}m_d ({E})^{1/\widetilde{p} }
\end{equation}
whenever $f$ is nonnegative, $\widetilde{p} >\frac{d}{k}$, and $E\subset
B(0,2)$. (The constant $C$ will depend on $\widetilde{p}$ and $\cC$.)
For an multi-index ${\bf n} =(n_1 ,\dots ,n_k )$ we will write
$\{|y_i | \sim 2^{n_i} \}$ to stand for the set of $y\in
\bbR ^k$ for which the $k$ inequalities $2^{n_i} \leq |y_i |<2^{n_i +1}$
hold. Our main task will be to establish the estimate
\begin{equation}\label{ineq6}\begin{aligned}
\int_{\bbR ^k} \int_{\{|y_i | \sim 1\}} f(x) \chi _E \big(y; L_1
(x,y),\dots ,L_l
(x,y)\big)\,
&
dm_k (y)\,dm_k (x)
\\
&
\leq C
\,\|f\|_{L^{d/k}(\bbR ^k )}m_d ({E})^{k/d}
\end{aligned}\end{equation}
for all nonnegative $f$ and $E\subset \bbR ^d$. From this a 
change of variables shows that the inequalities 
\begin{equation}\begin{aligned}\label{ineq7}
\int_{\bbR ^k} \int_{\{|y_i | \sim 2^{n_i} \}} f(x) \chi _E \big(y; L_1
(x,y),\dots
,L_l
(x,y)\big)\,
&
dm_k (y)\,dm_k
(x)
\\
&
\leq C
\,\|f\|_{L^{d/k}(\bbR ^k )}m_d ({E})^{k/d}
\end{aligned}\end{equation}
hold uniformly in $\bf n$. This implies \eqref{ineq5}: 
suppose $E\subset [-2, 2 ]^d$. For a multi-index $\bf n$ let 
$E_{\bf n}$ be the set of $(y_1 ,\dots ,y_k ;u_1 ,\dots ,u_l )\in E$ for
which $|y_i | \sim 2^{n_i}$. Then $m_d (E_{\bf n})\leq  
2^{\sum\nolimits_i (n_i +1)}$ and so, if
$\frac{1}{\widetilde{p}}=\frac{k}{d}-\epsilon$, 
$$
m_d (E_{\bf n})^{k/d}\leq m_d (E)^{1/\widetilde{p} 
}2^{\epsilon\sum\nolimits_i
(n_i +1)}. 
$$
Applying \eqref{ineq7} with $E$ replaced by $E_{\bf n}$ and then summing
over $\bf n$ for which $-\infty <n_j \leq 0$ yields \eqref{ineq5}.

Moving to the proof of \eqref{ineq6}, we write, for suitable functions $g$
on $\bbR ^ d$,  
\begin{equation}\begin{aligned}\label{defT}
\int_{\bbR ^k} \int_{\{|y_i | \sim 1\}} f(x) g \big(y; L_1
(x,y),
&
\dots 
,L_l (x,y)\big)\, 
 dm_k (y)\,dm_k (x)
=\langle Tf,g\rangle
\\
&
=
\int_{\bbR ^l} \int_{\{|y_i | \sim 1\}}
Tf(y;u)\,g(y;u)\,dm_k(y)\,dm_l (u).
\end{aligned}\end{equation}
Then \eqref{ineq6} is a consequence of the fact, which we will
establish below, that 
\begin{equation}
\label{Tbound}
T: L^{d/k}(\bbR ^k )\rightarrow L^{d/l}(\bbR ^d ).
\end{equation}
Although it does not figure here, one can regard the
operator $T$ as a restricted $(k-l)$-plane transform
operating on a function $f$ defined on $\bbR ^k$ by integrating $f$ over
the $(k-l)$-plane  
$$
\{x\in \bbR ^k :L_1 (x,y)=u_1 ,\dots ,L_l (x,y)=u_l\}.
$$
Since the indices in \eqref{Tbound} are conjugate, it is natural to
attempt to prove \eqref{Tbound} by embedding $T$ in an analytic family of
operators $\{T_z \}$ and then interpolating between $L^1 \rightarrow
L^{\infty}$ and $L^2 \rightarrow L^2$ estimates. Thus we define 
$$
T_z f(y;u)=C(z)f(y;\cdot\,)\ast |\cdot |^z (u),
$$
where the convolution is in the $u$ variable and $C(z)$ is chosen to
compensate for the singularities of the distributions $|\cdot |^z$ on 
$\bbR ^l$ --
see p. 363 in \cite{GS}. Next we will observe that
\begin{equation}\label{1inftybound}
\|T_z f\|_{L^{\infty}(\bbR ^d )}\leq c_0 (y)\, \|f\|_{L^{1}(\bbR ^k )}
\end{equation}
if $z=0+is$ and then prove (using the hypothesis (*) )that 
\begin{equation}\label{22bound}
\|T_z f\|_{L^{2}(\bbR ^d )}\leq c_1 (y)\, \|f\|_{L^{2}(\bbR ^k )}    
\end{equation}
if $z=-\frac{d}{2}+is$. 

Note that \eqref{defT} implies that $Tf(y;u)=0$ unless $|y_j |\sim 1$. If
$|y_j |\sim 1$ we will need the following formula:
\begin{equation}\label{formula}
\int_{\bbR ^l} Tf(y;u) h(u)\,dm_l (u)=\int_{\bbR ^k} f(x)\,h\big(L_1
(x,y),\dots ,L_l (x,y) 
\big)\, dm_k (x),
\end{equation}
valid for nice functions $h$ on $\bbR ^l$. To see
\eqref{formula} with $y=\widetilde{y}$, fix $\widetilde{y}$ with
$|\widetilde{y}_i |\sim 1$,
take $g(y;u)=\chi _{B(\widetilde{y},\delta )}(y)h(u)$ in the extreme terms
of 
\eqref{defT}, and then let $\delta \rightarrow 0$.  

Now \eqref{1inftybound} follows immediately from \eqref{formula}. 

To prove \eqref{22bound} we start by setting some notation. For fixed 
$y\in \bbR ^k$ we consider the mapping $\cL _y$ of $\bbR ^k$ into $\bbR
^l$ defined by 
$$
\cL_y x=\big( L_1 (x,y),\dots ,L_l (x,y)\big)
$$
along with the adjoint map $\cL^*_y$ of $\bbR ^l$ to $\bbR ^k$ defined by
$$
\langle \cL^*_y \zeta ,x\rangle =\langle \zeta ,\cL_y x \rangle .
$$
Then \eqref{formula} implies that 
\begin{equation}\label{FT}
\widehat{Tf (y,\cdot\,)}(\zeta )=\widehat{f}(\cL^*_y \zeta ).
\end{equation}
In order to prove \eqref{22bound} by exploiting \eqref{FT}, we need a
lemma.

\noindent{\bf Lemma.} {\it
Under the assumption (*) on $\cC$, there is
$c$, depending on $\cC$ and $\rho\in\bbR$, such that the inequality 
\begin{equation}\label{integralformula}
\int_{\{|y_j |\sim 1\}}\int_{\bbR ^l} |\zeta |^{\rho} w(\cL_y^* \zeta
)\,dm_l
(\zeta ) \,dm_k (y)
\leq c\, \int_{\bbR ^k} |\tau |^{\rho -k+l}w(\tau )\, dm_k (\tau )
\end{equation}
holds for nonnegative functions $w$ on $\bbR ^k$.
}

\noindent{\it Proof of Lemma:}
If $x,y\in\bbR ^k$, we may write
$x(i)$ instead of
$x_i$ and $xy$ to stand for the vector with $xy(i)=x(i)y(i)$. Let $\bf 1$
stand for
the vector $(1,1,\dots ,1)$. One may check that, for $i=1,\dots ,k$, 
$\cL_{\bf 1}^* \zeta (i)= \sum\nolimits_j c_i^j \zeta _j$ and also that 
$\cL_y^* \zeta =y\cL_{\bf 1}^* \zeta$. In particular, the hypothesis (*)
on $\cC$ has the following interpretation in terms of the $\binom k{l}$
coordinate projections $\pi$ of $\bbR ^k$ onto $\bbR ^{l}$:
for each such $\pi$, $\pi \circ \cL_{\bf 1}^* : \bbR ^l \rightarrow \bbR
^{l}$ is nonsingular. It follows that there is $M<\infty$ such that 
if $P\subset\{1,\dots ,k\}$ satisfies $|P|=l$, then 
\begin{equation}\label{normest}
|\zeta |\leq M \,\sup_{i\in P}|\cL_{\bf 1}^* \zeta (i)|.
\end{equation}
Next note that if $\zeta \in \bbR ^l$, then there are 
$1\leq i_{l+1}<i_{l+2}<\cdots <i_k\leq k$, dependent on $\zeta$, such that 
\begin{equation}
\label{normest2}
|\zeta |\leq M\, |\cL_{\bf 1}^* \zeta (i_a )|\  \text{if}\ a=l+1,\dots ,k.
\end{equation} 
(Having $|\zeta | > M\,|\cL_{\bf 1}^* \zeta (i)|$ for
even $l$ $i$'s would contradict \eqref{normest}.) In this situation,
write $Q=\{i_{l+1},\dots ,i_k \}$ and $\zeta \in F_Q$ so that 
$\bbR ^l =\cup_Q F_Q$, where the union is taken over all
$Q\subset\{1,\dots ,k\}$ such that $|Q|=k-l$. Then
\eqref{integralformula} will follow by summing over $Q$ the estimates
\begin{equation}\label{integralformula2}
\int_{\{|y_i |\sim 1\}}\int_{F_Q} |\zeta |^{\rho} w(\cL_y^* \zeta )\,dm_l
(\zeta
) \,
dm_k (y)
\leq c\, \int_{\bbR ^k} |\tau |^{\rho -k+l}w(\tau )\, dm_k (\tau ).
\end{equation}
To establish \eqref{integralformula2}, fix first $Q=\{i_{l+1},\dots
,i_k \}$, then 
$i_1 ,\dots ,i_{l}$ with $\{i_1 ,\dots ,i_k\}=\{1,\dots ,k\}$, and
finally $y_{i_1},\dots ,y_{i_{l}}$ with $|y_{i_a}|\sim 1$. 
Consider the map 
\begin{equation}\label{map}
(\zeta _1 ,\dots ,\zeta _l,y_{i_{l+1}},\dots ,y_{i_k})\mapsto \tau
\doteq \big(y_{i_1}\cL_{\bf 1}^*\zeta (i_1),\dots ,y_{i_k}\cL_{\bf
1}^*\zeta (i_k)\big)\backsimeq \cL_{y}^* \zeta ,
\end{equation}
where the $\backsimeq$ indicates a permutation of the coordinates. 
We want to estimate the absolute value $J$ of the Jacobian determinant
of \eqref{map} when $\zeta\in F_Q$. To do this, write $\tau$ as 
$$
\Big(y_{i_1}\sum c^j_{i_1}\zeta _j,\dots ,
y_{i_l}\sum c^j_{i_l}\zeta _j,
y_{i_{l+1}}\sum c^j_{i_{l+1}}\zeta _j,\dots,
y_{i_k}\sum c^j_{i_k}\zeta _j\Big),
$$
where $\sum$ means $\sum_{j=1}^l$. Computing the Jacobian matrix, one sees 
that
$$
J=\prod_{a=1}^l |y_{i_a}|\times |D(i_1 ,\dots ,i_l )| \times
\prod_{a=l+1}^k |\cL_{\bf 1}^* \zeta (i_a )|, 
$$
where $D(i_1 ,\dots ,i_l )$ is the determinant of the $l$ by $l$ matrix
obtained by retaining only the rows of $\cC$ corresponding to $i=i_1
,\dots ,i_l$. By (*), $|D(i_1 ,\dots ,i_l )|\geq c(\cC
)>0$. Since $\zeta \in F_Q$, 
$|\cL_{\bf 1}^* \zeta (i_a )|\geq \frac{|\zeta |}{M}$ for $a=l+1,\dots,
k$. It then follows from $|y_i |\sim 1$ that 
\begin{equation}\label{Jest}
J\geq c\, |\zeta |^{k-l}.
\end{equation}
It is also easy to check (see \eqref{normest}) that $|\zeta |^{\rho} \leq
c\,|\cL_{\bf 1}^* \zeta |^{\rho}$. So the inequality 
$$
\int_{\{|y_i |\sim 1\}}\int_{F_Q} |\zeta |^{\rho} \, w(\cL_y^* \zeta
)\,dm_l
(\zeta
) \,
dy_{i_{l+1}}\cdots dy_{i_k}
\leq c\, \int_{\bbR ^k} |\tau |^{\rho -k+l}w(\tau )\, dm_k (\tau )
$$
follows by change of variables, and then \eqref{integralformula2}
follows by integrating 
with respect to $y_{i_1},\dots ,y_{i_l}$ (since $|y_i |\sim 1$). 
This concludes the proof of the lemma.

With \eqref{integralformula} we can prove \eqref{22bound}: suppose $z=
-\frac{d}{2}+is$. Then 
\begin{equation*}\begin{aligned}
\int_{\{|y_i |\sim 1\}} \int_{\bbR ^l} |T_z f(y,u)|^2 
&
dm_l (u) 
\,dm_k (y)
=\int_{\{|y_i |\sim 1\}} \int_{\bbR ^l} |\widehat{T_z f(y,\cdot )}(\zeta
)|^2
dm_l (\zeta) \,dm_k (y)
\\
&
=c(s)\int_{\{|y_i |\sim 1\}}\int_{\bbR ^l}|\widehat{f}(\cL_y^* \zeta )|^2
\big| |\zeta |^{d/2-l}\big| ^2 dm_l (\zeta)\, dm_k (y)
\\
&
\leq c(s)\int_{\bbR ^k}|\widehat{f}(\tau )|^2 |\tau |^{d-2l-(k-l)}dm_k
(\tau )
=c(s)\|f\|_{L^2 (\bbR ^k )}^2 ,
\end{aligned}\end{equation*}
where: the second equality follows from \eqref{FT} and the
fact that, on $\bbR ^l$, $\widehat{|\cdot |^z}(\zeta )=c(z)\,|\zeta
|^{-z-l}$ (\cite{GS}, p. 363); the inequality follows from
\eqref{integralformula}; and the
last equality follows from $d=k+l$. This proves \eqref{22bound}. Now
interpolating between \eqref{1inftybound} and \eqref{22bound} shows that 
\begin{equation*}
T_z :L^{d/k}(\bbR ^k )\rightarrow L^{d/l}(\bbR ^d )
\end{equation*} 
if $z=-l+is$. Since $T_{-l}$ is a scalar multiple of $T$, \eqref{Tbound}
follows, concluding the proof of the theorem.

\end{document}